\documentclass[11pt]{article}
\usepackage{amsmath, amsthm,amsfonts,amssymb,amscd}
\usepackage{color}
\usepackage{geometry}
\usepackage{tikz}
\usepackage{subfigure}
\numberwithin{equation}{section}
\usepackage[T1]{fontenc}
\usepackage{cite}
\usepackage{enumerate}

\usepackage{float}
\usepackage{soul}
\usepackage{marginnote}
\newcommand{\Limsup}{\mathop{{\rm Lim}\,{\rm sup}}}

\def\tto{\;{\lower 1pt \hbox{$\rightarrow$}}\kern -10pt
\hbox{\raise 2pt \hbox{$\rightarrow$}}\;}

\def\ra{\rangle}
\def\la{\langle}

\def\h{\hfill\Box}
\def\R{\mathbb R}

\def\ox{\bar{x}}

\def\oy{\bar{y}}

\def\h{\hfill\triangle}

\newcounter{lk}
\def\Limsup{\mathop{{\rm Lim}\,{\rm sup}}}

\begin{document}

\newtheorem{Theorem}{Theorem}[section]
\newtheorem{Proposition}[Theorem]{Proposition}
\newtheorem{Remark}[Theorem]{Remark}
\newtheorem{Lemma}[Theorem]{Lemma}
\newtheorem{Corollary}[Theorem]{Corollary}
\newtheorem{Definition}[Theorem]{Definition}
\newtheorem{Example}[Theorem]{Example}
\newtheorem{Fact}[Theorem]{Fact}
\renewcommand{\theequation}{\thesection.\arabic{equation}}
\normalsize
\def\proof{
\normalfont
\medskip
{\noindent\itshape Proof.\hspace*{6pt}\ignorespaces}}
\def\endproof{$\h$ \vspace*{0.1in}}

\title{\small \bf ON THE CODERIVATIVE OF THE PROJECTION OPERATOR\\ ONTO THE POSITIVE CONE IN HILBERT SPACES}
\date{}

\author{Le Van Hien\footnote{Faculty of Pedagogy, HaTinh University, Hatinh, Vietnam; email: hien.levan@htu.edu.vn.}\ \ \ and  \ Nguyen Viet Quan\footnote{Faculty of Pedagogy, HaTinh University, Hatinh, Vietnam; email: nguyenvietquan231097@gmail.com.}}

\maketitle

{\small \begin{abstract}

In this paper, we study the generalized differentiability of the metric projection operator onto the positive cone in Hilbert spaces. We first establish the formula for exactly computing the regular coderivative and the Mordukhovich coderivative of the metric projection operator onto the positive cone in Euclidean spaces.  Then, these results are also established for the projection operator onto the positive cone in the real Hilbert space $l_2$.

\end{abstract}}

\medskip

\normalsize
\section{Introduction}
\setcounter{equation}{0}

Let  $H$ be a Hilbert space with norm $\| . \|$
and let $C$ be a nonempty closed convex subset of $H$. The metric projection operator $P_C$ is defined by
$$P_C(x)=\{y\in C:\inf\limits_{z\in C}\|x-z\|=\|x-y\|\},\ \ \mbox{for all}\ \ x\in H.$$

The metric projection operator $P_C$ is a well-defined single-valued mapping that has many useful properties, such as continuity, monotonicity, and non-expansiveness. Researching the metric projection operator in Hilbert space is an interesting topic and receives the attention of many people, because of its applications in optimization, theory of equation, control theory and many other fields; see \cite{FP82, H77, KL24, K84, Li23.3, Li24.3, Li23, Li23.2, Li24, Li24.2, Li24.4, LLX23, M02, N95, S87, WY03} and the references therein.

In addition to the above properties, the smoothness of the projection operator is also of special interest. Recently, in their articles, Li  has studied the Gâteaux directional differentiability of the projection operator in Bochner space \cite{Li23.3} and the projection operator in the real uniformly convex and uniformly smooth Banach space \cite{Li23}, in this paper, Li gave the exact representations of the directional derivatives of $P_C$ with $C$ is a closed ball or a closed and convex cone (including proper closed subspaces); studied the strict Fréchet differentiability of the projection operator in Hilbert spaces \cite{Li24}, the Fréchet differentiability of the projection operator in uniformly convex and uniformly smooth Banach spaces \cite{Li24.3}. In that direction of research, Hien \cite{H24} also studied the strict Fréchet differentiability of the projection operator onto closed balls with centerat arbitrarily given point in Hilbert spaces and onto the second-order cones in Euclidean spaces.

Generalized differentiation lies at the heart of variational analysis andits applications. In \cite{M18}, Mordukhovich introduced the concept of generalized differentiation and provide many useful properties, of which the coderivative is an important tool because of its applications in optimization and algorithms. Therefore, research on the coderivative is one of the interesting topics, attracting the attention of many people; see \cite{DR14, HOS12, LM04, M94, M18, MO07, OS08, RW98, YZ17} and the references therein. Recently, calculating the coderivative of the metric projection operator has also obtained some interesting results. In \cite{OS08}, Outrata and Sun  have calculated the coderivative of the metric projection onto the second-order cone in $\R^n$, then, this result was used to obtain a sufficientcondition for the Aubin property of the solution map of a parameterized secondorder cone complementarity problem and to derive necessary optimality conditionsfor a mathematical program with a second-order cone complementarity problemamong the constraints. Li \cite{Li23.2, Li24.2, Li24.4}  investigated the properties of the regular coderivative of the metric projection operator onto some closed balls, closed and convex cylinders and positive cones in Hilbert spaces, in uniformly convex and uniformly smooth Banach spaces, and in some general Banach spaces. In \cite{H24(2)}, Hien gave the complete calculation formula for the regular coderivative of the projection operator onto closed balls with center at arbitrarily given point in Hilbert spaces.

The aim of this paper is to establish a formula for calculating the coderivative  of the metric projection operator onto the positive cone in Hilebert spaces. This problem has been studied by Li \cite{Li23.2}, however the results on the regular coderivative of the metric operator projected onto the positive cones in \cite{Li23.2} is only given in some special cases and has an inaccurate conclusion. More precisely, we first establish aformula for exactly computing the regular coderivative and Mordukhovich coderivative  of the metric projection operator onto the positive cone in Euclidean space $\R^n$. Then, these formulas are also established for the projection operator onto the positive cone in the real Hilbert space $l_2$.

The structure of the paper is as follows.   In section 2,  we recall some preliminary materials. In section 3, we present results on the regular coderivative and Mordukhovich coderivative of the metric projection operator onto the positive cone in  $\R^n$. Then, in section 4, we give the exact calculation formula for the regular coderivative and Mordukhovich coderivative of the metric projection operator onto the positive cone in $l_2$. Finally, we conclude the paper in section 5 where we discuss some perspectives of the obtained results and future works.

\section{Preliminaries}
\setcounter{equation}{0}
We first give the following notations that will be used throughout the paper. Let $X$ be Hilbert space  with scalar product $\la . , . \ra$ and norm $\|.\|$. Given a set $C \subset X,$ we denote by $C^o$ its interior.  Denote $t\downarrow 0$ means that $t\to 0$ with $t>0$. The symbol $\mathbb{N}$ is the set of natural numbers.

Below are basic notions and facts from variational analysis, which are frequently used in the following; see \cite{A99, Li24, M18, RW98, S96} for more details.

First, We recall the concepts of Gâteaux directional differentiability and the Fréchet differentiability of the metric projection operator $P_C:H\rightarrow C.$ For $\ox\in X$ and $w\in X\backslash\{\theta\}$, if the following limit exists,
$$P_C^{'}(\ox)(w)=\lim\limits_{t\downarrow0}\dfrac{P_C(\ox+tw)-P_C(\ox)}{t},$$ then, $P_C$ is said to be Gâteaux directionally differentiable at point $\ox$ along direction $w$ and $P_C^{'}(\ox)(w)$ is called the Gâteaux directional derivative of $P_C$ at point $\ox$ along direction $w$. 

$P_C$ is called Fréchet differentiable at $\ox \in X$ if there is a linear continuous mapping $\nabla P_C(\ox):X\to X$ such that
$$\lim\limits_{h\to 0}\dfrac{\|P_C(\ox+h)-P_C(\ox)-\nabla P_C(\ox)(h)\|}{\|h\|}=0.$$
Then, $\nabla P_C(\ox)$  is said to be the  Fréchet	derivative of $P_C$ at $\ox$.

In particular, the mapping $P_C$ is called  strictly Fréchet differentiable at $\ox \in X$ if
$$\lim\limits_{(u,v)\to (\ox,\ox)}\dfrac{\|P_C(u)-P_C(v)-\nabla P_C(\ox)(u-v)\|}{\|u-v\|}=0.$$

The mapping $P_C$ is called  Gâteaux directionally differentiable (Fréchet differentiable, strictly Fréchet differentiable) on $\Omega\subset X$ if $P_C$ is Gâteaux directionally differentiable (Fréchet differentiable, strictly Fréchet differentiable, respectively) at every point $x\in \Omega$.

\begin{Definition} {\rm (see \cite{M18,RW98})
		Let $\Omega$ be a nonempty subset of  a Hilbert space $H$ and $\ox\in\Omega$, we denote by $\widehat N_\Omega(\ox)$ the {\it Fréchet (regular) normal cone} to $\Omega$ at $\ox$, defined by
		$$\widehat N_\Omega(\ox):=\Big\{z\in H| \limsup\limits_{u\xrightarrow{\Omega} \ox}\dfrac{\la z,u-\ox\ra}{\|u-\ox\|} \leq 0 \Big\}.$$
The {\it limiting (Mordukhovich) normal cone} to $\Omega$ at $\ox$, denoted $N_\Omega(\ox)$, is defined by 
$$N_\Omega(\ox):=\Limsup\limits_{u\xrightarrow{\Omega} \ox}\widehat N_\Omega(x)$$
where "Lim sup" is the Painlevé-Kuratowski outer limit of sets (see \cite{RW98})	
}	
\end{Definition}
If $\ox\not\in\Omega$, one puts
$\widehat N_\Omega(\ox)=N_\Omega(\ox)=\emptyset$ by convention. When the set $\Omega$ is convex, $\widehat N_\Omega(\ox)=N_\Omega(\ox)$ amounts to the classic normal cone in the sense of
convex analysis.

\begin{Definition} {\rm (see \cite{M18,RW98})
		Let $H$ be Hilbert space and $\Phi: H\rightrightarrows H$ be a set-valued mapping with its graph {\rm gph}$\Phi: =\{(x,y)|y\in\Phi(x)\}$ and its domain {\rm Dom}$\Phi:=\{x|\Phi(x)\not=\emptyset\}$. The multifunction $\widehat D^\ast\Phi(\ox,\oy):H\rightrightarrows H$, defined by
			$$\widehat D^\ast\Phi(\ox,\oy)(w):=\big\{v\in H| (v,-w)\in \widehat N_{\mbox{\rm gph}\Phi}(\ox, \oy)\}\ \mbox{for all}\ w\in H$$ 
	is called {\it regular coderivative} of $\Phi$ at $(\ox, \oy).$
		
Analogously, the multifunction
$D^\ast\Phi(\ox,\oy):H\rightrightarrows H$, defined by
$$D^\ast\Phi(\ox,\oy)(w):=\big\{v\in H| (v,-w)\in  N_{\mbox{\rm gph}\Phi}(\ox, \oy)\}\ \mbox{for all}\ w\in H$$ 
is called {\it limiting (Mordukhovich) coderivative} of $\Phi$  at $(\ox, \oy).$ 		
	 }
\end{Definition}
In the case $\Phi(\ox)=\{\oy\},$ we usually write   $\widehat D^\ast\Phi(\ox) (D^\ast\Phi(\ox)).$ If $\Phi: H\rightarrow H$ is strictly Fréchet differentiable at $\ox$, then $D^\ast\Phi(\ox)(y)=\widehat D^\ast\Phi(\ox)(y)=\{\nabla\Phi(\ox)(y)\}, \ \mbox{for all}\ y\in H.$\\

Note that {\rm gph}$\Phi \subset H\times H$, where $H\times H$ be the orthogonal product of $H$ equipped with the inner product $\la. , .\ra_{H\times H}$, which is abbreviated as $\la. , .\ra$.  It is defined by, $$\la(x,y), (u,v)\ra=\la x,v\ra +\la y,v\ra,\ \ \mbox{for any}\  (x,y), (u,v)\in H \times H.$$
Let $\| . \|_{H\times H}$ be the norm in $H\times H$ induced by the inner product $\la. , .\ra$. Then, we have
$$\| (x,y) \|_{H\times H}=\sqrt{\la(x,y), (x,y)\ra} =\sqrt{\|x\|^2+\|y\|^2}\ \ \mbox{for all}\ x,y\in H.$$

The following result is given in \cite{H24(2)}, which provides a useful formula to calculate the Fréchet coderivatives of single-valued mappings on Hilbert spaces.
 \begin{Lemma} {\rm (see \cite{H24(2)})}\label{Lem1}
 	Let $H$ be a Hilbert space, $f: H\rightarrow H$ be a Lipschitz  continuous mapping on $H$. Then, the regular coderivatives of $f$ at a point $\ox\in H$ satisfies that, for any $y\in H$,
 	$$\widehat D^\ast f(\ox)(y)=\left\{z\in H\big|\limsup\limits_{u\rightarrow \ox}\dfrac{\la z,u-\ox\ra-\la y,f(u)-f(\ox)\ra}{\|u-\ox\|} \leq 0 \right\}.$$
 \end{Lemma}

 \section{Coderivative of the metric projection onto the positive cone in $\R^n$}\label{Sec3}
\setcounter{equation}{0}
Let $\mathbb{R}^n$ denote the n-dimensional Euclidean space with origin $\theta$. We recall some notations and concepts  from \cite{Li24}, which will be used in this section.  For any $x=(x_1,...,x_n)\in\R^n$, we define three subsets of the set $\{1,2,...,n\}$ with respect to the given $x$ by
$$x^+:=\{i\in \{1,...,n\}: x_i>0\};$$
$$x^-:=\{i\in \{1,...,n\}: x_i<0\};$$
$$x^\bullet:=\{i\in \{1,...,n\}: x_i=0\}.$$
Let K denote the positive cone of $\R^n$, which is defined by
$$K:=\{x\in\R^n:x^-=\emptyset\}.$$
Then, $K$ is a pointed closed and convex cone in $\R^n$ with the interior given by
$$K^o=\{x\in\R^n:x^-=x^\bullet=\emptyset\}.$$
Define $$\widehat K:= \{x\in\R^n: x^+\not=\emptyset, x^-\not=\emptyset, x^\bullet=\emptyset\},$$
and $$\Delta\R^n=\{x\in\R^n: x^\bullet\not=\emptyset\}.$$
For any given fixed $x\in\widehat K$, we define a mapping $b(x;.):\R^n\rightarrow \R^n$, by
$$b(x;w)_i:=\begin{cases}
w_i\ \ \ \ \mbox{if}\ i\in x^+,\\
0 \ \ \ \ \ \ \mbox{if}\ i\not\in x^+,
\end{cases} \ \ \mbox{for all}\ w\in \R^n.$$
It follows from \cite{Li24} that  the metric projection operator $P_K:\R^n\rightarrow K$ has the following properties:
 
(i) $$P_K(x)_i:=\begin{cases}
	x_i\ \ \ \ \mbox{if}\ i\in x^+,\\
	0 \ \ \  \ \ \mbox{if}\ i\not\in x^+,
\end{cases} \ \ \mbox{for all}\ x\in \R^n.$$

(ii) $P_K$ is a Lipschitz continuous mapping on $\R^n$ with Lipschitz constant $L=1$, i.e
 $$\|P_K(x)-P_K(y)\|\leq \|x-y\|,\ \ \mbox{for all} \ x,y\in\R^n.$$

(iii) $P_K$ is strictly Fréchet differentiable on $K^o, (-K)^o$ and $\widehat K$; $P_{K}$ is Gâteaux directionally differentiable but not Fréchet differentiable on $\Delta\R^n$.

The following result provides the formula for calculating the regular coderivatives of the metric projection onto the positive cone in $\R^n$.

	\begin{Theorem} \label{Thm1} 
	Let $K$ be the positive cone of $\R^n$. Then, the regular coderivative of the projection operator $P_K$ at $\ox\in\R^n$ is given by
	\begin{equation}\label{kq}\widehat D^\ast P_K(\ox)(y)=\left\{z\in \R^n: z_i=y_i, i\in \ox^+; z_i=0, i\in \ox^-; 0\leq z_i\leq y_i, i\in \ox^\bullet\right\}, \ \ \mbox{for all} \ y\in\R^n.\end{equation}

	\end{Theorem}
	
	{\bf Proof.}
	
  Let $y\in \R^n$, thanks to the Lipschitz continuity of $P_{K}$, by Lemma \ref{Lem1}, we have
\begin{equation}\label{39}z\in \widehat D^\ast P_{K}(\ox)(y) \iff \limsup\limits_{u\rightarrow \ox}\dfrac{\la z,u-\ox\ra-\la y,P_{K}(u)-P_{K}(\ox)\ra}{\|u-\ox\|} \leq 0.\end{equation}

Suppose that $z\in \widehat D^\ast P_{K}(\ox)(y)$. We take a directional line segment in the limit in \eqref{39}, $u =(u_i)\in\R^n$  with $u_i=\ox_i+
t(z_i-y_i), \ \mbox{for}\ i\in \ox^+; \ u_i=\ox_i+
tz_i,\ \mbox{for}\ i\in \ox^-;\ u_i=0, \ \mbox{for} \ i\in \ox^\bullet$, where $t\downarrow 0$. It follows that
$$\begin{array}{rl}
	0&\geq \limsup\limits_{u\rightarrow \ox}\dfrac{\la z,u-\ox\ra-\la y,P_{K}(u)-P_{K}(\ox)\ra}{\|u-\ox\|}\\
	&\geq \limsup\limits_{t\downarrow 0}\dfrac{\sum\limits_{i\in\ox^+}(z_i-y_i)(u_i-\ox_i)+\sum\limits_{i\in\ox^-}z_i(u_i-\ox_i)+\sum\limits_{i\in\ox^\bullet}(z_iu_i-y_i(P_K(u))_i)}{\|u-\ox\|}\\
		&= \limsup\limits_{t\downarrow 0}\dfrac{t\sum\limits_{i\in\ox^+}(z_i-y_i)^2+t\sum\limits_{i\in\ox^-}z_i^2}{\sqrt{t^2\sum\limits_{i\in\ox^+}(z_i-y_i)^2+t^2\sum\limits_{i\in\ox^-}z_i^2}}\\
	&=\sqrt{\sum\limits_{i\in\ox^+}(z_i-y_i)^2+\sum\limits_{i\in\ox^-}z_i^2}.
	\end{array}$$
This implies that \begin{equation}\label{y1} z_i=y_i,\ \mbox{for all}\ i\in \ox^+ \ \ \mbox{and}\ \ z_i=0,\ \mbox{for all}\ i\in \ox^-.\end{equation} 

Assume, by the way of contradiction, that that there exists $i_0 \in\ox^\bullet$ such that $z_{i_0}>y_{i_0}$. Choosing $u =(u_i)\in\R^n$  with $u_i=\ox_i, \ \mbox{for }\ i\not= i_0; \ u_{i_0}=t$ for $t\downarrow 0$. We have
$$\begin{array}{rl}
	0&\geq \limsup\limits_{u\rightarrow \ox}\dfrac{\la z,u-\ox\ra-\la y,P_{K}(u)-P_{K}(\ox)\ra}{\|u-\ox\|}\\
	&\geq \limsup\limits_{t\downarrow 0}\dfrac{\sum\limits_{i\in\ox^+}(z_i-y_i)(u_i-\ox_i)+\sum\limits_{i\in\ox^-}z_i(u_i-\ox_i)+\sum\limits_{i\in\ox^\bullet}(z_iu_i-y_i(P_K(u))_i)}{\|u-\ox\|}\\
	&= \limsup\limits_{t\downarrow 0}\dfrac{t\left(z_{i_0}-y_{i_0}\right)}{t}\\
	&=z_{i_0}-y_{i_0}.
\end{array}$$
This means $z_{i_0}\leq y_{i_0}$. This contradiction proves that 
\begin{equation} \label{y2}
z_i\leq y_i \  \ \ \ \ \mbox{for all}\ i\in \ox^\bullet.
\end{equation}

Similarly, suppose there exists $i_0 \in\ox^\bullet$ such that $z_{i_0}<0$. By choosing $u =(u_i)\in\R^n$  with $u_i=\ox_i, \ \mbox{for }\ i\not= i_0; \ u_{i_0}=-t$ for $t\downarrow 0$. We have
$$\begin{array}{rl}
	0&\geq \limsup\limits_{u\rightarrow \ox}\dfrac{\la z,u-\ox\ra-\la y,P_{K}(u)-P_{K}(\ox)\ra}{\|u-\ox\|}\\
	&\geq \limsup\limits_{t\downarrow 0}\dfrac{\sum\limits_{i\in\ox^+}(z_i-y_i)(u_i-\ox_i)+\sum\limits_{i\in\ox^-}z_i(u_i-\ox_i)+\sum\limits_{i\in\ox^\ast}(z_iu_i-y_i(P_K(u))_i)}{\|u-\ox\|}\\
	&= \limsup\limits_{t\downarrow 0}\dfrac{-tz_{i_0}}{t}\\
	&=-z_{i_0}.
\end{array}$$
This shows that \begin{equation} \label{y3}
z_i\geq 0, \ \ \ \mbox{for all} \ i\in \ox^\bullet. 
\end{equation}

By \eqref{y1}, \eqref{y2} and \eqref{y3} together, we have 
\begin{equation} \label{c}
	z\in \widehat D^\ast P_{K}(\ox)(y)\Longrightarrow  \begin{cases} z_i=y_i\ \ \ \ \ \ \ \ \mbox{if}\ \ i\in \ox^+\\
		 z_i=0\ \ \ \ \ \ \ \ \  \mbox{if} \ i\in \ox^-\\
		  0\leq z_i\leq y_i \ \ \ \mbox{if}\  i\in \ox^\bullet 
		  \end{cases}
\end{equation}

Next, we prove the converse, assuming $z=(z_i)\in\R^n$ with $z_i=y_i,  \ i\in \ox^+;$$ 
z_i=0,\ i\in\ox^- $ and $0\leq z_i\leq y_i, \  i\in \ox^\bullet$. We need to show that
$$\limsup\limits_{u\rightarrow \ox}\dfrac{\la z,u-\ox\ra-\la y,P_{K}(u)-P_{K}(\ox)\ra}{\|u-\ox\|}\leq 0.$$

Indeed, we have 
$$\begin{array} {rl}&\dfrac{\la z,u-\ox\ra-\la y,P_{K}(u)-P_{K}(\ox)\ra}{\|u-\ox\|}\\
	&=\dfrac{\sum\limits_{i\in\ox^+}(z_i-y_i)(u_i-\ox_i)+\sum\limits_{i\in\ox^-}z_i(u_i-\ox_i)+\sum\limits_{i\in\ox^\bullet}(z_iu_i-y_i(P_K(u))_i)}{\|u-\ox\|}\\
	&=\dfrac{\sum\limits_{i\in\ox^\bullet}(z_iu_i-y_i(P_K(u))_i)}{\|u-\ox\|}\\
	&=\dfrac{\sum\limits_{i\in\ox^\bullet\cap u^+}(z_i-y_i)u_i+\sum\limits_{i\in\ox^\bullet\backslash u^+}z_iu_i}{\|u-\ox\|}\\
	&\leq 0 \ \ \ \ \mbox{for all}\ u\in\R^n\backslash\{\ox\}.
		\end{array}$$

So, $\limsup\limits_{u\rightarrow \ox}\dfrac{\la z,u-\ox\ra-\la y,P_{K}(u)-P_{K}(\ox)\ra}{\|u-\ox\|}\leq 0$. This means \begin{equation} \label{d}
\begin{cases} z_i=y_i\ \ \ \ \ \ \ \ \mbox{if}\ \ i\in \ox^+\\
	z_i=0\ \ \ \ \ \ \ \ \  \mbox{if} \ i\in \ox^-\\
	0\leq z_i\leq y_i \ \ \ \mbox{if}\  i\in \ox^\bullet 
\end{cases}\Longrightarrow 	z\in \widehat D^\ast P_{K}(\ox)(y).
\end{equation}
Then, by \eqref{c} and \eqref{d}, \eqref{kq} is proved.\hfill $\square$
\\

Using Theorem \ref{Thm1}, we obtain the following result, part of which was given by  Li in \cite{Li23.2}.
\begin{Corollary} \label{Co1} {\rm (see \cite[Theorem 4.1]{Li23.2}) }
Let $K$ be the positive cone of $\R^n$ with negative cone $-K$. Then, the regular coderivative of the projection operator $P_K$ have the following representations. 

(i) For any $y\in\R^n$, we have$$\widehat D^\ast P_K(\ox)(y)=\begin{cases}
	\{y\}\ \ \ \ \ \ \ \  \ \mbox{if}\ \ox\in K^o\\
	\{\theta\}\ \ \ \ \ \ \ \ \ \mbox{if}\ \ox\in (-K)^o\\
	\{b(\ox;y)\}\ \ \ \mbox{if}\ \ox\in \widehat K\\
\end{cases}$$

(ii) For any $\ox\in \Delta\R^n$, we have

(a) $ \widehat D^\ast P_{K}(\ox)(\theta)=\{\theta\};$\\

(b) For any $y\in \R^n$, 
$$y^-\cap\ox^\bullet\not=\emptyset\Longrightarrow \lambda y\not\in \widehat D^\ast P_{K}(\ox)(y)  \ \mbox{for all}\ \lambda<1.$$

In particular, 
$$y^-\cap\ox^\bullet\not=\emptyset\Longrightarrow \theta\not\in \widehat D^\ast P_{K}(\ox)(y).$$

(c) For any $y\in \R^n$, $y^-\cap\ox^\bullet\not=\emptyset\Longrightarrow \widehat D^\ast P_{K}(\ox)(y) =\emptyset.$

(d)  $\widehat D^\ast P_{K}(\ox)(\ox)=\{P_K(\ox)\}$.

(e)  $\widehat D^\ast P_{K}(\theta)(y)=\displaystyle\prod_{i=1}^n[0, y_i].$
 
\end{Corollary}
	{\bf Proof.}
	(i) For any $\ox\in K^o$, we have $\ox^-=\ox^\bullet=\emptyset$. For any $y\in \R^n$, using \eqref{kq}, we get $$\begin{array}{rl} \widehat D^\ast P_K(\ox)(y)&=\left\{z\in \R^n: z_i=y_i, i\in \ox^+\right\}.\\
	&=\{z\in \R^n: z_i=y_i, \forall i=1,2,...,n\}\\
	&=\{y\}.
	\end{array}$$
	
For any $\ox\in (-K)^o$, we have $\ox^+=\ox^\bullet=\emptyset$. For any $y\in \R^n$, using \eqref{kq}, we get $$\begin{array}{rl} \widehat D^\ast P_K(\ox)(y)&=\left\{z\in \R^n: z_i=0, i\in \ox^-\right\}.\\
	&=\{z\in \R^n: z_i=0, \forall i=1,2,...,n\}\\
	&=\{\theta\}.
\end{array}$$

For any $\ox\in \widehat K$, we have $\ox^\bullet=\emptyset$. For any $y\in \R^n$, using \eqref{kq}, we get $$\begin{array}{rl} \widehat D^\ast P_K(\ox)(y)&=\left\{z\in \R^n: z_i=y_i, i\in \ox^+; z_i=0, i\in \ox^-\right\}.\\
	&=\left\{z\in \R^n: z_i=y_i, i\in \ox^+; z_i=0, i\not\in \ox^+\right\}.\\
	&=\{b(\ox;y)\}.
\end{array}$$	
	
(ii.a) Using \eqref{kq} with $y=0$, we have  
$$\begin{array}{rl}
\widehat D^\ast P_K(\ox)(\theta)&=\left\{z\in \R^n: z_i=0, i\in \ox^+; z_i=0, i\in \ox^-; 0\leq z_i\leq 0, i\in \ox^\bullet\right\}\\
&=\{z\in \R^n: z_i=0, \forall i=1,2,...,n\}\\
&=\{\theta\}.
\end{array}$$

(ii.b)  If $y^-\cap\ox^\bullet\not=\emptyset$, then there exists $i_0\in \ox^\bullet$ such that $y_{i_0}<0$. We have $\lambda y_{i_0}>y_{i_0}$ for all $\lambda<1$. By \eqref{kq}, $\lambda y\not\in \widehat D^\ast P_{K}(\ox)(y)  \ \mbox{for all}\ \lambda<1.$ In particular, for $\lambda =0$, we get $\theta\not\in \widehat D^\ast P_{K}(\ox)(y).$

(ii.c) Since $y^-\cap\ox^\bullet\not=\emptyset$, there exists $i_0\in \ox^\bullet$ such that $y_{i_0}<0$. This means that there does not exist  $z=(z_i)\in\R^n$ such that $0\leq z_{i_0}\leq y_{i_0}$, and $\widehat D^\ast P_{K}(\ox)(y) =\emptyset.$

(ii.d) Using \eqref{kq} with $y=\ox$, we have  
$$\begin{array}{rl}
	\widehat D^\ast P_K(\ox)(\ox)&=\left\{z\in \R^n: z_i=\ox_i, i\in \ox^+; z_i=0, i\in \ox^-; 0\leq z_i\leq 0, i\in \ox^\bullet\right\}\\
	&=\left\{z\in \R^n: z_i=\ox_i, i\in \ox^+; z_i=0, i\not\in \ox^+\right\}\\
	&=\{P_K(\ox)\}.
\end{array}$$

(ii.e) Using \eqref{kq} with $\ox=\theta$, we have  $\ox^+=\ox^-=\emptyset, \ \ox^\bullet =\{1,2,...,n\}$ and
$$\begin{array}{rl}
	\widehat D^\ast P_K(\theta)(y)&=\left\{z\in \R^n:  0\leq z_i\leq y_i, i\in \ox^\bullet\right\}.\\
	&=\left\{z\in \R^n: 0\leq z_i\leq y_i, \forall i=1,2,...,n\right\}\\
	&=\displaystyle\prod_{i=1}^n[0, y_i].
\end{array}$$\hfill $\square$

\begin{Remark}
By Corollary \ref{Co1} (ii.d), the statements of \cite[Theorem 4.1 (iv.c)]{Li23.2} 

"Let $\ox\in\Delta\R^n\backslash\{\theta\}$. For any $z\in\R^n$, we have,
\begin{equation}
	z\not=\theta \Longrightarrow z\not\in \widehat D^\ast P_{K}(\ox)(\ox),
\end{equation}
\begin{equation}
	\mbox{and}\ \ 	\ox^+\not=\emptyset \Longrightarrow \widehat D^\ast P_{K}(\ox)(\ox)=\emptyset.
\end{equation}"

are incorrect.
\end{Remark}

Next, we compute the limiting (Mordukhovich) coderivative $D^\ast P_{K}(\ox)$. Since $P_K(.)$ is continuous, the graph of $P_K$ is closed and, by \cite[Equation 8(18)]{RW98}, we know that 
\begin{equation}\label{dn}D^\ast P_{K}(\ox)(y)=\Limsup\limits_{x\rightarrow \ox, y'\rightarrow y}\widehat D^\ast P_{K}(x)(y')\end{equation}

This, together with  Theorem \ref{Thm1}, allows us to provide a complete characterization of $D^\ast P_{K}(\ox)$.
\begin{Theorem} \label{Thm34}
	Let $K$ be the positive cone of $\R^n$. Then, the Mordukhovich coderivative of the projection operator $P_K$ at $\ox\in\R^n$ is given by
\begin{equation}\label{kq2}\begin{array}{rl}D^\ast P_K(\ox)(y)&=\widehat D^\ast P_K(\ox)(y)\\
		&=\left\{z\in \R^n: z_i=y_i, i\in \ox^+; z_i=0, i\in \ox^-; 0\leq z_i\leq y_i, i\in \ox^\bullet\right\},\end{array}
	\end{equation} $\mbox{for all} \ y\in\R^n.$
\end{Theorem}

{\bf Proof.} Since $\widehat D^\ast P_K(\ox)(y)\subset D^\ast P_K(\ox)(y)$, it suffices to prove
$$D^\ast P_K(\ox)(y)\subset \widehat D^\ast P_K(\ox)(y).$$
Take $z\in D^\ast P_K(\ox)(y)$. Then by definition, there exist $y^k\rightarrow y, x^k\rightarrow \ox$ and $z^k\rightarrow z$ such that $z^k\in \widehat D^\ast P_K(x^k)(y^k)$, for all $k\in\mathbb N.$ By \eqref{kq}, we have $$\begin{cases} z^k_i=y^k_i\ \ \ \ \ \ \ \ \mbox{if}\ \ i\in (x^k)^+\\
z^k_i=0\ \ \ \ \ \ \ \ \ \ \mbox{if} \ i\in (x^k)^-\\
0\leq z^k_i\leq y^k_i \ \ \ \mbox{if}\  i\in (x^k)^\bullet 
\end{cases}\ \ \ \ \mbox{for all}\ k\in\mathbb{ N}.$$
Since $x^k\rightarrow \ox$, we have $\ox^+\subset (x^k)^+, \ox^-\subset (x^k)^-$ and $(x^k)^\bullet\subset \ox^\bullet$. So, for all $i\in \ox^+\subset (x^k)^+$, we have  \begin{equation} \label{1} z_i=\lim\limits_{k\rightarrow\infty}z^k_i=\lim\limits_{k\rightarrow\infty}y^k_i=y_i.\end{equation} 
For all $i\in \ox^-\subset (x^k)^-$, \begin{equation} \label{2} z_i=\lim\limits_{k\rightarrow\infty}z^k_i=\lim\limits_{k\rightarrow\infty}0=0.\end{equation}
 For $i\in \ox^\bullet$, using a subsequence if necessary, we consider the following cases.

{\it Case 1.} $i\in (x^k)^+,$ for all $k\in\mathbb{N}$. We have $z_i=\lim\limits_{k\rightarrow\infty}z^k_i=\lim\limits_{k\rightarrow\infty}y^k_i=y_i.$

{\it Case 2.} $i\in (x^k)^-,$ for all $k\in\mathbb{N}$. We have $z_i=\lim\limits_{k\rightarrow\infty}z^k_i=\lim\limits_{k\rightarrow\infty}0=0.$

{\it Case 3.} $i\in (x^k)^\bullet,$ for all $k\in\mathbb{N}$. We have $0\leq z^k_i \leq y^k_i$ for all $k\in \mathbb{N}.$ Taking the limits as $k\rightarrow\infty$, we get $0\leq z_i \leq y_i$.

Thus, for any $i\in \ox^\bullet$, \begin{equation}\label{3}0\leq z_i \leq y_i\end{equation}
Then, by \eqref{1}, \eqref{2}, \eqref{3} and \eqref{kq}, $z\in \widehat D^\ast P_K(\ox)(y).$ The proof is complete. \hfill $\square$\\

To end this section, we give some examples to demonstrate the results of Theorems \ref{Thm1} and \ref{Thm34}.
\begin{Example} {\rm 
	Let $\ox=(\ox_1,\ox_2)\in\R^2$. For all $y=(y_1,y_2)\in \R^2,$ we have
	$$\widehat D^\ast P_K(\ox)(y)=D^\ast P_K(\ox)(y)=\begin{cases}
		\{y\}\ \ \ \ \ \ \ \ \ \ \ \ \ \ \ \ \  \ \ \mbox{if}\ \ \ox_1, \ox_2 >0,\\
		\left\{(0, 0)\right\}  \ \ \ \ \ \ \ \ \ \ \  \ \ \mbox{if}\ \ \ox_1, \ox_2 <0,\\
			\left\{(\ox_1, 0)\right\} \ \ \ \ \ \ \ \ \ \ \ \  \mbox{if}\ \ \ox_1>0, \ox_2<0,\\
				\left\{(0, \ox_2)\right\} \ \ \ \ \ \ \ \ \ \ \ \  \mbox{if}\ \ \ox_1<0, \ox_2>0,\\
		\left\{y_1\right\}\times [0, y_2]\  \ \ \ \ \ \  \mbox{if}\ \ \ox_1>0, \ox_2=0,\\
 [0, y_1]\times	\left\{y_2\right\}\  \ \ \ \ \ \  \mbox{if}\ \ \ox_1=0, \ox_2>0,\\
		\left\{0\right\}\times [0, y_2]\ \ \ \ \ \ \ \ \  \mbox{if}\ \ \ox_1<0, \ox_2=0,\\
		[0, y_1]\times	\left\{0\right\}\ \ \ \ \ \ \ \ \  \mbox{if}\ \ \ox_1=0, \ox_2<0,\\
		[0, y_1]\times	[0,y_2]\  \ \ \ \ \ \  \mbox{if}\ \ \ox_1=\ox_2=0.
	\end{cases}$$
}
\end{Example}

\begin{Example} {\rm
		Let  $a,b,c>0.$ Then, for any $y=(y_1,y_2,y_3)\in\R^3$ we have
	$$\widehat D^\ast P_K(a,b,c)(y)=D^\ast P_K(a,b,c)(y)=\{y\}.$$
		$$\widehat D^\ast P_K(-a,-b,-c)(y)=D^\ast P_K(-a,-b,-c)(y)=\{(0,0,0)\}.$$
			$$\widehat D^\ast P_K(a,b,0)(y)=D^\ast P_K(a,b,0)(y)=\{y_1\}\times \{y_2\}\times [0,y_3].$$
			$$\widehat D^\ast P_K(a,b,-c)(y)=D^\ast P_K(a,b,-c)(y)=\{(y_1, y_2, 0)\}.$$	
			$$\widehat D^\ast P_K(a,-b,-c)(y)=D^\ast P_K(a,-b,-c)(y)=\{(y_1, 0, 0)\}.$$
			$$\widehat D^\ast P_K(0,-b,-c)(y)=D^\ast P_K(0,-b,-c)(y)=[0, y_1]\times \{0\}\times \{0\}.$$	
				$$\widehat D^\ast P_K(a,0,-c)(y)=D^\ast P_K(a,0,-c)(y)=\{y_1\}\times [0,y_2]\times \{0\}.$$
					$$\widehat D^\ast P_K(a,0,0)(y)=D^\ast P_K(a,0,0)(y)=\{y_1\}\times [0, y_2]\times [0,y_3].$$	
				$$\widehat D^\ast P_K(0,0,-c)(y)=D^\ast P_K(0,0,-c)(y)=[0,y_1]\times [0,y_2]\times \{0\}.$$
				$$\widehat D^\ast P_K(0,0,0)(y)=D^\ast P_K(0,0,0)(y)=[0,y_1]\times [0,y_2]\times [0,y_3].$$ }
\end{Example}

\section{Coderivative of the metric projection onto the positive cone in real Hilbert space $l_2$}\label{Sec4}
\setcounter{equation}{0}

In this section, we consider the real Hilbert space $l_2$ with norm $\| . \|$ and inner product $\la .,. \ra$. Let $\mathbb K$ denote the positive cone of $l_2$ defined by
$$\mathbb{K}:=\left\{x=(x_1,x_2,...)\in l_2: x_i\geq 0, \ \mbox{for all}\ i\in \mathbb{N}\right\}.$$
Let $P_\mathbb{K}: l_2\rightarrow \mathbb{K}$ be the metric projection operator from $l_2$ onto $\mathbb{K}$. It follows from \cite{Li24} that for any $x\in l_2$, $P_\mathbb{K}(x)$ is represented as follows
$$P_\mathbb{K}(x)_i=\begin{cases}
x_i\ \ \ \ \ \mbox{if}\ x_i>0;\\
0\ \ \ \ \ \ \mbox{if}\ x_i\leq0;
\end{cases} \mbox{for}\  i\in\mathbb N.$$
Furthermore, $P_\mathbb{K}$ is a Lipschitz continuous mapping on $l_2$ with Lipschitz constant $L=1$, i.e
$$\|P_\mathbb{K}(x)-P_\mathbb{K}(y)\|\leq \|x-y\|, \ \mbox{for all }\ x,y\in l_2.$$
And $P_\mathbb{K}$ is positive homogeneous, i.e
$$P_\mathbb{K}(\lambda x)=\lambda P_\mathbb{K}(x),\ \ \mbox{for any}\ \ x\in l_2\ \mbox{and}\ \lambda>0.$$

By \cite[Theorem 5.3]{Li24}, $P_\mathbb{K}$ is Gâteaux directionally differentiable but not Fréchet differentiable on sets $\mathbb{K}^+, \mathbb{K}^-$ and $\widehat{\mathbb{K}}$. Where, 
$$\mathbb{K}^+:=\left\{x=(x_1,x_2,...)\in l_2: x_i> 0, \ \mbox{for all}\ i\in \mathbb{N}\right\};$$
$$\mathbb{K}^-:=\left\{x=(x_1,x_2,...)\in l_2: x_i< 0, \ \mbox{for all}\ i\in \mathbb{N}\right\};$$
$$\begin{array}{rl}\widehat{\mathbb{K}}:=\{&x=(x_1,x_2,...)\in l_2: |x_i|> 0, \ \mbox{for all}\ i\in \mathbb{N}\  \mbox{and}\\ &
		 \mbox{there are at least one pair} \ j, k\in\mathbb{N}\ \mbox{with}\ x_jx_k<0  \}.\end{array}$$
For any $x=(x_1,x_2,...)\in l_2$, define
$$x^+:=\{i\in \mathbb{N}: x_i>0\};$$
$$x^-:=\{i\in \mathbb{N}: x_i<0\};$$
$$x^\bullet:=\{i\in \mathbb{N}: x_i<0\}.$$
The following result provides the formula for calculating the regular coderivatives of the metric projection onto the positive cone in $l_2$.

\begin{Theorem} \label{Thm41} 
	Let $\mathbb K$ be the positive cone of $l_2$. Then, the regular coderivative of the projection operator $P_\mathbb K$ at $\ox\in l_2$ is given by
	\begin{equation}\label{kq22}\widehat D^\ast P_\mathbb K(\ox)(y)=\left\{z\in l_2: z_i=y_i, i\in \ox^+; z_i=0, i\in \ox^-; 0\leq z_i\leq y_i, i\in \ox^\bullet\right\}, \ \ \mbox{for all} \ y\in l_2.\end{equation}
	
\end{Theorem}

{\bf Proof.} Similar to the proof of Theorem \ref{Thm1}. Take $\ox\in l_2$.

Let $y\in l_2$, thanks to the Lipschitz continuity of $P_{K}$, by Lemma \ref{Lem1}, we have
\begin{equation}\label{30}z\in \widehat D^\ast P_{K}(\ox)(y) \iff \limsup\limits_{u\rightarrow \ox}\dfrac{\la z,u-\ox\ra-\la y,P_{K}(u)-P_{K}(\ox)\ra}{\|u-\ox\|} \leq 0.\end{equation}

Suppose that $z\in \widehat D^\ast P_{K}(\ox)(y)$, we prove 
$$\begin{cases} z_i=y_i\ \ \ \ \ \ \ \ \mbox{if}\ \ i\in \ox^+,\\
	z_i=0\ \ \ \ \ \ \ \ \  \mbox{if} \ i\in \ox^-,\\
	0\leq z_i\leq y_i \ \ \ \mbox{if}\  i\in \ox^\bullet.
\end{cases}$$
Indeed, suppose there exists $j\in \ox^+$ such that $z_j\not=y_j,$ by choosing $u\in l_2$ with $u_i=\ox_i,$ for $i\not=j$ and $u_j=x_j+t(z_j-y_j),$ for $t\downarrow 0.$ We have $u \rightarrow \ox$ as $t\downarrow 0$ and 
$$\begin{array}{rl}
	0&\geq \limsup\limits_{u\rightarrow \ox}\dfrac{\la z,u-\ox\ra-\la y,P_{K}(u)-P_{K}(\ox)\ra}{\|u-\ox\|}\\
	&\geq \limsup\limits_{t\downarrow 0}\dfrac{\sum\limits_{i\in\ox^+}(z_i-y_i)(u_i-\ox_i)+\sum\limits_{i\in\ox^-}z_i(u_i-\ox_i)+\sum\limits_{i\in\ox^\bullet}(z_iu_i-y_i(P_K(u))_i)}{\|u-\ox\|}\\
	&= \limsup\limits_{t\downarrow 0}\dfrac{t(z_j-y_j)^2}{t|z_j-y_j|}\\
	&=|z_j-y_j|.
\end{array}$$
This implies that \begin{equation}\label{v1} z_i=y_i,\ \mbox{for all}\ i\in \ox^+.\end{equation} 

Suppose there exists $j\in \ox^-$ such that $z_j\not=0,$ by choosing $u\in l_2$ with $u_i=\ox_i,$ for $i\not=j$ and $u_j=x_j+tz_j,$ for $t\downarrow 0.$ We have $u \rightarrow \ox$ as $t\downarrow 0$ and 
$$\begin{array}{rl}
	0&\geq \limsup\limits_{u\rightarrow \ox}\dfrac{\la z,u-\ox\ra-\la y,P_{K}(u)-P_{K}(\ox)\ra}{\|u-\ox\|}\\
	&\geq \limsup\limits_{t\downarrow 0}\dfrac{\sum\limits_{i\in\ox^+}(z_i-y_i)(u_i-\ox_i)+\sum\limits_{i\in\ox^-}z_i(u_i-\ox_i)+\sum\limits_{i\in\ox^\bullet}(z_iu_i-y_i(P_K(u))_i)}{\|u-\ox\|}\\
	&= \limsup\limits_{t\downarrow 0}\dfrac{tz_j^2}{t|z_j|}\\
	&=|z_j|.
\end{array}$$
This implies that \begin{equation}\label{v2} z_i=0,\ \mbox{for all}\ i\in \ox^-.\end{equation} 
Assume, by the way of contradiction, that that there exists $j \in\ox^\bullet$ such that $z_j>y_j$. Choosing $u =(u_i)\in l_2$  with $u_i=\ox_i$, for $i\not= j$ and  $u_j=t$, for $t\downarrow 0$. We have $u\rightarrow \ox$ and
$$\begin{array}{rl}
	0&\geq \limsup\limits_{u\rightarrow \ox}\dfrac{\la z,u-\ox\ra-\la y,P_{K}(u)-P_{K}(\ox)\ra}{\|u-\ox\|}\\
	&\geq \limsup\limits_{t\downarrow 0}\dfrac{\sum\limits_{i\in\ox^+}(z_i-y_i)(u_i-\ox_i)+\sum\limits_{i\in\ox^-}z_i(u_i-\ox_i)+\sum\limits_{i\in\ox^\ast}(z_iu_i-y_i(P_K(u))_i)}{\|u-\ox\|}\\
	&\geq \limsup\limits_{t\downarrow 0}\dfrac{(z_ju_j-y_j(P_K(u))_j)}{t}\\
	&= \limsup\limits_{t\downarrow 0}\dfrac{t\left(z_j-y_j\right)}{t}\\
	&=z_j-y_j.
\end{array}$$
This means $z_j\leq y_j$. This contradiction proves that 
\begin{equation} \label{v3}
	z_i\leq y_i, \  \ \ \ \ \mbox{for all}\ i\in \ox^\bullet.
\end{equation}
Similarly, suppose there exists $j \in\ox^\bullet$ such that $z_j<0$. By choosing $u =(u_i)\in l_2$  with $u_i=\ox_i,$ for $j\not=j$ and $u_j=-t$, for $t\downarrow 0.$ We get $u\rightarrow \ox$ and
$$\begin{array}{rl}
	0&\geq \limsup\limits_{u\rightarrow \ox}\dfrac{\la z,u-\ox\ra-\la y,P_{K}(u)-P_{K}(\ox)\ra}{\|u-\ox\|}\\
	&\geq \limsup\limits_{t\downarrow 0}\dfrac{\sum\limits_{i\in\ox^+}(z_i-y_i)(u_i-\ox_i)+\sum\limits_{i\in\ox^-}z_i(u_i-\ox_i)+\sum\limits_{i\in\ox^\ast}(z_iu_i-y_i(P_K(u))_i)}{\|u-\ox\|}\\
	&\geq \limsup\limits_{t\downarrow 0}\dfrac{(z_ju_j-y_j(P_K(u))_j)}{t}\\
	&= \limsup\limits_{t\downarrow 0}\dfrac{-tz_j}{t}\\
	&=-z_j.
\end{array}$$
This shows that \begin{equation} \label{v4}
	z_i\geq 0, \ \ \ \mbox{for all} \ i\in \ox^\bullet. 
\end{equation}

By \eqref{v1}, \eqref{v2}, \eqref{v3} and \eqref{v4} together, we have 
\begin{equation} \label{c2}
	z\in \widehat D^\ast P_{K}(\ox)(y)\Longrightarrow  \begin{cases} z_i=y_i\ \ \ \ \ \ \ \ \mbox{if}\ \ i\in \ox^+\\
		z_i=0\ \ \ \ \ \ \ \ \  \mbox{if} \ i\in \ox^-\\
		0\leq z_i\leq y_i \ \ \ \mbox{if}\  i\in \ox^\bullet 
	\end{cases}
\end{equation}

Next, we prove the converse, assuming $z\in l_2$ with $z_i=y_i$, for $i\in \ox^+,$ $z_i=0,$ for $i\in \ox^-$ and $0\leq z_i\leq y_i,$ for  $i\in~\ox^\bullet$. We need to show that
$$\limsup\limits_{u\rightarrow \ox}\dfrac{\la z,u-\ox\ra-\la y,P_{K}(u)-P_{K}(\ox)\ra}{\|u-\ox\|}\leq 0.$$

Indeed, we have 
$$\begin{array} {rl}&\dfrac{\la z,u-\ox\ra-\la y,P_{K}(u)-P_{K}(\ox)\ra}{\|u-\ox\|}\\
	&=\dfrac{\sum\limits_{i\in\ox^+}(z_i-y_i)(u_i-\ox_i)+\sum\limits_{i\in\ox^-}z_i(u_i-\ox_i)+\sum\limits_{i\in\ox^\bullet}(z_iu_i-y_i(P_K(u))_i)}{\|u-\ox\|}\\
	&=\dfrac{\sum\limits_{i\in\ox^\bullet}(z_iu_i-y_i(P_K(u))_i)}{\|u-\ox\|}\\
	&=\dfrac{\sum\limits_{i\in\ox^\bullet\cap u^+}(z_i-y_i)u_i+\sum\limits_{i\in\ox^\bullet\backslash u^+}z_iu_i}{\|u-\ox\|}\\
	&\leq 0 \ \ \ \ \mbox{for all}\ u\in l_2\backslash\{\ox\}.
\end{array}$$

So, $\limsup\limits_{u\rightarrow \ox}\dfrac{\la z,u-\ox\ra-\la y,P_{K}(u)-P_{K}(\ox)\ra}{\|u-\ox\|}\leq 0$. This means \begin{equation} \label{d2}
	\begin{cases} z_i=y_i\ \ \ \ \ \ \ \ \mbox{if}\ \ i\in \ox^+\\
		z_i=0\ \ \ \ \ \ \ \ \  \mbox{if} \ i\in \ox^-\\
		0\leq z_i\leq y_i \ \ \ \mbox{if}\  i\in \ox^\bullet 
	\end{cases}\Longrightarrow 	z\in \widehat D^\ast P_{K}(\ox)(y).
\end{equation}
Then, by \eqref{c2} and \eqref{d2}, \eqref{kq22} is proved.\hfill $\square$\\

Let $N$ be a nonempty subset of $\mathbb{N}$ with complementary $\bar N$. We define some subsets and  an ordering relation in $l_2$ with respect to $N$.
$$\mathbb{R}^N=\left\{x=(x_i)\in l_2: x_i=0, \ \mbox{for}\ i\in \bar N\right\};$$
$$\mathbb{K}_N=\left\{x=(x_i)\in l_2: x_i\geq0, \ \mbox{for}\ i\in  N\right\};$$
$$\partial\mathbb{K}_N=\left\{x=(x_i)\in \mathbb{K}_N: x_i= 0, \ \mbox{for}\ i\in  N\right\};$$
$$\mathbb{Z}_N=\left\{x=(x_i)\in l_2: x_i> 0, \ \mbox{for all}\ i\in  N \ \mbox{and}\ x_i=0, \ \mbox{for all}\ i\in  \bar N \right\};$$
$$z\preceq_N y \Leftrightarrow z_i\leq y_i, \ \mbox{for all}\ i\in  N \ \mbox{and}\ z_i=y_i, \ \mbox{for all}\ i\in  \bar N.$$

Using Theorem \ref{Thm41}, we obtain the following result, which was given by  Li in \cite{Li23.2}.

\begin{Corollary} {\rm (see \cite[Theorem 5.1]{Li23.2}) }
Let $M$ be a nonempty finite subset of $\mathbb{N}$ with complementary $\bar M$. The regular coderivative of $P_{\mathbb{K}}$ have the following representations.
	
	(i) For any $\ox\in l_2$, we have $$\widehat D^\ast P_\mathbb K(\ox)(\theta)=\{\theta\};$$
	
	(ii) Let $\ox\in \mathbb{Z}_M.$ For any $y\in l_2$, 
\begin{equation} \label{hq1}
y\in \mathbb{K}_{\bar M} \Longleftrightarrow y\in \widehat D^\ast P_\mathbb K(\ox)(y);
\end{equation}	

(iii) Let $\ox\in \mathbb{Z}_M.$ For any $y\in \mathbb{K}_{\bar M}$, we have
\begin{equation} \label{hq2}
\widehat D^\ast P_\mathbb K(\ox)(y) =\left\{z\in \mathbb{K}_{\bar M}: z\preceq_M y \right\}.
\end{equation}	
In particular,
$$y\in \partial\mathbb{K}_{\bar M} \Longrightarrow \widehat D^\ast P_\mathbb K(\ox)(y) =\{y\}.$$	
\end{Corollary}
{\bf Proof.}

(i)  Using \eqref{kq22} with $y=\theta$, we have  
$$\begin{array}{rl}
	\widehat D^\ast P_\mathbb K(\ox)(\theta)&=\left\{z\in l_2: z_i=0, i\in \ox^+; z_i=0, i\in \ox^-; 0\leq z_i\leq 0, i\in \ox^\bullet\right\}\\
	&=\{z\in l_2: z_i=0, \forall i\in\mathbb{N}\}\\
	&=\{\theta\}.
\end{array}$$

(ii, iii)  Since $\ox\in \mathbb{Z}_M$, $\ox^+=M, \ox^-=\emptyset$ and $\ox^\bullet=\bar M.$ By \eqref{kq22}, for any $y\in l_2$, we have
$$\begin{array}{rl}
	\widehat D^\ast P_\mathbb K(\ox)(y)&=\left\{z\in l_2: z_i=y_i, i\in \ox^+; z_i=0, i\in \ox^-; 0\leq z_i\leq y_i, i\in \ox^\bullet\right\}\\
	&=\left\{z\in l_2: z_i=y_i, i\in M;  0\leq z_i\leq y_i, i\in \bar M\right\}\\
	&=\{z\in l_2: z\preceq_{\bar M} y,\ z_i\geq 0, i\in\bar M\}\\
	&=\left\{z\in \mathbb{K}_{\bar M}: z\preceq_{\bar M} y \right\}.
\end{array}$$
So, $$\begin{array} {rl} y\in \widehat D^\ast P_\mathbb K(\ox)(y) &\Longleftrightarrow y\in \left\{z\in \mathbb{K}_{\bar M}: z\preceq_{\bar M} y \right\}\\
	& \Longleftrightarrow y\in \mathbb{K}_{\bar M}.\end{array}$$
	In particular, 
	$$\begin{array} {rl}y\in \partial\mathbb{K}_{\bar M} &\Longrightarrow y_i=0, \ \mbox{for all}\ i\in\bar M\\
		&\Longrightarrow \widehat D^\ast P_\mathbb K(\ox)(y)=\left\{z\in \mathbb{K}_{\bar M}: z\preceq_{\bar M} y \right\}\\
		&\Longrightarrow \widehat D^\ast P_\mathbb K(\ox)(y)=\{y\}.
		\end{array}
	$$
	\hfill $\square$\\
	
Similar to Theorem \ref{Thm34}, the following result is deduced from Theorem \ref{Thm41} and the definition of the Mordukhovich coderivative.
	
\begin{Theorem} \label{Thm43}
	Let $\mathbb K$ be the positive cone of $l_2$. Then, the Mordukhovich coderivative of the projection operator $P_\mathbb K$ at $\ox\in l_2$ is given by
	\begin{equation}\label{kq2}\begin{array}{rl}D^\ast P_\mathbb K(\ox)(y)&=\widehat D^\ast P_\mathbb K(\ox)(y)\\
			&=\left\{z\in l_2: z_i=y_i, i\in \ox^+; z_i=0, i\in \ox^-; 0\leq z_i\leq y_i, i\in \ox^\ast\right\},\end{array}
	\end{equation} $\mbox{for all} \ y\in l_2.$
\end{Theorem}	
{\bf Proof.}	 Proof of this Theorem is similar to the proof of Theorem \ref{Thm34}. Therefore, it is omitted here.

\section{Concluding remarks}

The main results of this paper exhibits formulas for computing the regular coderivative and Mordukhovich coderivative of the metric projection operator onto positive cones in Euclidean spaces and in real Hilbert space $l_2$. In the near future, we will continue to research other generalized differential structures of metric projection operators in Hilbert spaces and general Banach spaces.

.
\section*{Acknowledgments}
The first author would like to thank Vietnam Institute for Advanced Study in Mathematics for hospitality during his post-doctoral fellowship of the Institute in 2022--2023.


\small
\begin{thebibliography}{100}
\bibitem{A99}  R.   Akerkar,  Nonlinear Functional Analysis, American Mathematical Society. 1999.
 


\bibitem{DR14} A. L. Dontchev and R. T. Rockafellar, Implicit Functions and Solution Mappings. A View from Variational Analysis, 2nd ed., Springer Ser. Oper. Res. Financ. Eng., Springer, NewYork, 2014.


 	
 \bibitem{FP82} S. Fitzpatrick and R. R. Phelps,  Differentiability of the metric projection in Hilbert space, Trans. Amer. Math.
Soc., 270(1982), 483-501.

\bibitem{H77} A. Haraux,  How to differentiate the projection on a convex set in Hilbert space, Some applications to variational inequalities, J. Math. Soc. Japan., 29(1977) 615-631.

\bibitem{HOS12} R. Henrion, J. V. Outrata and T. Surowiec,  On regular coderivatives in parametric equilibria with non-unique multipliers, Mathematical Programming., 136(2012), 111-131.

\bibitem{H24(2)}  L. V. Hien, Regular coderivative and graphical derivative of the metric projection onto closed balls in Hilbert spaces, 2024. https://arxiv.org/abs/2406.18377


\bibitem{H24}  L. V. Hien, Some results on the strict Fréchet differentiability of the metric projection operator in Hilbert spaces, 2024. https://arxiv.org/abs/2403.14512


\bibitem{KL24} A. A. Khan  and J. L. Li,  Characterizations of the Metric and Generalized Metric Projections on Subspaces of Banach Spaces, J. Math. Anal. Appl., 531(2024).


\bibitem{K84} A. Korányi,  Monotone functions on formally real Jordan algebras. Math. Ann., 269(1984), 73-76. 

\bibitem{LM04} A.B. Levy and B.S. Mordukhovich, Coderivatives in parametric optimization. Math. Program. 99(2004), 311-327.

\bibitem{Li23.3} J. L. Li, \newblock Directional Differentiability of the Metric Projection in Bochner Spaces, (2023). https://arxiv.org/abs/2311.00942

\bibitem{Li24.3} J. L. Li, \newblock Fréchet differentiability of the metric projection operator in Banach spaces, (2024). https://arxiv.org/abs/2401.01480v1

\bibitem{Li23} J. L. Li, \newblock Directional differentiability of the metric projection operator in uniformly convex and uniformly smooth Banach spaces, to appear in \newblock{ Journal of Optimization Theory and Applications}, (2024).

\bibitem{Li24} J. L. Li, \newblock Strict Fréchet differentiability of the metric projection operator in Hilbert spaces, (2024). https://arxiv.org/abs/2312.14362

\bibitem{Li23.2} J. L. Li, \newblock Mordukhovich derivatives of the metric projection operator in Hilbert spaces, (2023).
https://arxiv.org/pdf/2401.01906

\bibitem{Li24.2} J. L. Li, \newblock Mordukhovich derivatives of the metric projection operator in uniformly convex and uniformly smooth Banach spaces, (2024). https://arxiv.org/abs/2401.11321

\bibitem{Li24.4} J. L. Li, \newblock Mordukhovich derivatives of the set-valued metric projection operator in general Banach spaces, (2024). https://arxiv.org/abs/2402.10827

\bibitem{LLX23} J. L. Li, L. Cheng, L. S. Liu and L. S. Xie, \newblock Directional differentiability of the metric projection in Hilbert spaces and Hilbertian Bochner spaces, to appear in \newblock{Journal of Nonlinear and Convex 	Analysis.}, (2024).

\bibitem{M02} K. Malanowski, Differentiability of projections onto cones and sensitivity analysis for optimal control,
Proceedings of the 41st IEEE Conference on Decision and Control, Las Vegas, Nevada USA, 2002.

\bibitem{M94} B. S. Mordukhovich,  \newblock Generalized differential calculus for nonsmooth and set-valued mappings,
J. Math. Anal. Appl. 183(1994), 250-288. 


\bibitem{M18} B. S. Mordukhovich,  \newblock Variational Analysis and Generalized Differentiation I, Basic Theory, Springer Heidelberg New York Dordrecht London, 2006.


\bibitem{MO07} B. S. Mordukhovich and J. V. Outrata, Coderivative analysis of quasivariational inequalities with applications to stability and optimization. SIAM J. Optim. 18(2007), 389-412.

\bibitem{N95} D. Noll, Directional differentiability of the metric projection in Hilbert space. Pacific Journal of Mathematics,  170(1995),  567-592. 


\bibitem{OS08} J. V. Outrata and D. Sun, \newblock On the Coderivative of the Projection Operator onto the Second-order Cone, \newblock{ Set-Valued Anal}, 16(2008), 999-1014. 

\bibitem{RW98} R. T. Rockafellar and   R. J-B. Wets,
\newblock { Variational analysis}. \newblock Springer-Verlag, Berlin, 1998.

\bibitem{S96}  E. Schechter, \newblock { Handbook of Analysis and its Foundations}. \newblock Academic Press, New York, 1997.

\bibitem{S87}   A. Shapiro, On differentiability of the metric projection in W1. Boundary case. Proc. Amer. Math. Soc.,  99(1987), 123-128.

\bibitem{WY03} Z. Wu, J. J. Ye, Equivalence among various derivatives and subdifferentials of the distance function, J. Math. Anal. Appl., 282(2003), 629-647.



\bibitem{YZ17} J. J. Ye, J. Zhou, Exact formulas for the proximal/regular/limiting normal cone of the second-order cone complementarity set, Mathematical Programming, Mathematical Programming, 162(2017), 33-50.






















\end{thebibliography}
\end{document}